\documentclass{article}
\usepackage{amsmath,amsthm,amsopn,amstext,amscd,amsfonts,amssymb}
\usepackage{natbib}
\usepackage{verbatim}

\def \diag {\mathop{\rm diag}\nolimits}
\def \build#1#2#3{\mathrel{\mathop{#1}\limits^{#2}_{#3}}}

\def \Prod#1#2{\mbox{$\displaystyle \build{\prod}{#2}{#1}$}}

\renewenvironment{abstract}
                 {\vspace{6pt}
                  \begin{center}
                  \begin{minipage}{5in}
                  \centerline{\textbf{Abstract}}
                  \noindent\ignorespaces
                 }
                 {\end{minipage}\end{center}}
\newtheorem{thm}{\textbf{Theorem}}[section]

\theoremstyle{definition}

\newtheorem{rem}{\textbf{Remark}}[section]

\title{\huge \textbf{More about measures and Jacobians of singular random matrices}}
\author{
  \textbf{Jos\'e A. D\'{\i}az-Garc\'{\i}a} \thanks{Corresponding author\newline
    }\\
  Department of Statistics and Computation \\
  25350 Buenavista, Saltillo, Coahuila, Mexico \\
  E-mail: jadiaz@uaaan.mx \\
}
\date{}
\begin{document}
\maketitle
\begin{abstract}
This work studies the Jacobians of certain singular transformations and the
corresponding measures which support the jacobian computations.
\end{abstract}

\section{Introduction}

First consider the following notation: Let ${\mathcal L}_{m,N}(q)$ be the linear
space of all $N \times m$ real matrices of rank $q \leq \min(N,m)$ and ${\mathcal
L}_{m,N}^{+}(q)$ be the linear space of all $N \times m$ real matrices of rank $q
\leq \min(N,m)$, with $q$ distinct singular values. The set of matrices
$\mathbf{H}_{1} \in {\mathcal L}_{m,N}$ such that $\mathbf{H}'_{1}\mathbf{H}_{1} =
\mathbf{I}_{m}$ is a manifold denoted ${\mathcal V}_{m,N}$, called Stiefel manifold.
In particular, ${\mathcal V}_{m,m}$ is the group of orthogonal matrices ${\mathcal
O}(m)$. Denote by ${\mathcal S}_{m}$, the homogeneous space of $m \times m$ positive
definite symmetric matrices; and by ${\mathcal S}_{m}^{+}(q)$, the $(mq - q(q -
1)/2)$-dimensional manifold of rank $q$ positive semidefinite $m \times m$ symmetric
matrices with $q$ distinct positive eigenvalues.

Assuming that $\mathbf{X} \in {\mathcal L}_{m,N}^{+}(q)$, \citet{dgm:97} proposed the
Jacobian of non-singular part of the singular value decomposition, $\mathbf{X} =
\mathbf{H}_{1}\mathbf{DW}'_{1}$, where $\mathbf{H}_{1} \in {\mathcal V}_{q,N}$,
$\mathbf{D}$ is a diagonal matrix with $D_{1}
> D_{2} > \cdots D_{q} > 0$ and $\mathbf{W}_{1} \in {\mathcal V}_{q,m}$. Also, note that the
jacobian itself defines the factorization of Hausdorff's measure $(d\mathbf{X})$ (or
Lebesgue's measure defined on the manifold ${\mathcal L}_{m,N}^{+}(q)$, see \citet[p.
249]{b:86}). Analogous results for $\mathbf{V} \in {\mathcal S}_{m}^{+}(q)$
considering the non-singular part of the spectral decomposition of $\mathbf{V}$ were
proposed by \citet{u:94} and \citet{dg:97}. Based on these two results,
\citet{dgt:05} and \citet{dgg:06} computed the jacobians of the transformations
$\mathbf{Y} = \mathbf{X}^{+}$ and $\mathbf{W} = \mathbf{V}^{+}$, where
$\mathbf{A}^{+}$ denotes the Moore-Penrose inverse of $A$, see \citet[p.49]{r:73}.

In the present work, assuming that $\mathbf{X} \in {\mathcal L}_{m,N}^{+}(q)$, we
proposed the Jacobian of non-singular part of the singular value decomposition
assuming multiplicity in the singular values of $\mathbf{X}$ and the corresponding
Jacobian of $\mathbf{Y} = \mathbf{X}^{+}$ under the same conditions. Analogous
results for $\mathbf{V} \in {\mathcal S}_{m}^{+}(q)$ and $\mathbf{W} =
\mathbf{V}^{+}$ considering the non-singular part of the spectral decomposition of
$\mathbf{V}$ are proposed assuming multiplicity in the eigenvalues of $\mathbf{V}$
and/or assuming that $\mathbf{V}$ is a indefinite singular matrix. Also we will
determine the explicit measures with respect the jacobians are computed, see
\citet{dgb:07}.

\section{Jacobian of symmetric matrices}\label{indef}

Consider again $\mathbf{A} \in \mathcal{S}_{m}$, it remains to study: $\mathbf{A}$ as
a (nonsingular) indefinite matrix, i.e. $\mathbf{A} \in \mathcal{S}_{m}^{\pm}(m_{1},
m_{2})$, with $m_{1} + m_{2} = m$, where $m_{1}$ is the number of positive
eigenvalues and $m_{2}$ is the number of negative eigenvalues; and $\mathbf{A}$ as a
(singular) semi-indefinite matrix, i.e. $\mathbf{A} \in \mathcal{S}_{m}^{\pm}(q,
q_{1}, q_{2})$, with $q_{1} + q_{2} = q$, here $q_{1}$ is the number of positive
eigenvalues and $q_{2}$ is the number of negative eigenvalues.

First suppose  $\mathbf{A} \in \mathcal{S}_{m}^{\pm}(m_{1}, m_{2})$ such that
$\mathbf{A} = \mathbf{HDH}'$, where $\mathbf{H} \in \mathcal{O}(m)$, $\mathbf{D}$ is
a diagonal matrix. Without loss of generality, let $\lambda_{1}>\cdots>
\lambda_{m_{1}}> 0$ and $0 > -\delta_{1}>\cdots> -\delta_{m_{2}}$, explicitly
$$
  \mathbf{A} = \mathbf{H} \left [
            \begin{array}{cccccc}
              \lambda_{1} & \cdots & 0 & 0 & \cdots & 0 \\
              \vdots & \ddots & \vdots & \vdots & \ddots & \vdots \\
              0 & \cdots & \lambda_{m_{1}} & 0 & \cdots & 0 \\
              0 & \cdots & 0 & -\delta_{1} & \cdots & 0 \\
              \vdots & \ddots & \vdots & \vdots & \ddots & \vdots \\
              0 & \cdots & 0 & 0 & \cdots & -\delta_{m_{2}}
            \end{array}
        \right ] \mathbf{H}'.
$$
Now let $\mathbf{A} \in \mathcal{S}_{3}^{\pm}(1,2)$ and let $A = \mathbf{HDH}'$ be
its SD, then
$$
  d\mathbf{A} = d\mathbf{HDH}' + \mathbf{H}d\mathbf{DH}' + \mathbf{HD}d\mathbf{H}',
$$
thus by the skew symmetry of $\mathbf{H}'d\mathbf{H}$ we have, see \citet[p.
105]{m:82}
$$
  \mathbf{H}'d\mathbf{AH} = \mathbf{H}'d\mathbf{HD} + d\mathbf{D} + \mathbf{D}d\mathbf{H}'\mathbf{H}
    = \mathbf{H}'d\mathbf{HD} + d\mathbf{D} - \mathbf{DH}'d\mathbf{H}
$$
Moreover,
\begin{eqnarray*}
  \mathbf{H}'d\mathbf{AH} &=& \left [
             \begin{array}{ccc}
               0 & -h'_{2}dh_{1} & -h'_{3}dh_{1} \\
               h'_{2}dh_{1} & 0 & -h'_{3}dh_{2} \\
               h'_{3}dh_{1} & h'_{3}dh_{2} & 0
             \end{array}
          \right]
          \left [
             \begin{array}{ccc}
               \lambda_{1} & 0 & 0 \\
               0 & -\delta_{1} & 0 \\
               0 & 0 & -\delta_{2}
             \end{array}
          \right] +
          \left [
             \begin{array}{ccc}
               d\lambda_{1} & 0 & 0 \\
               0 & -d\delta_{1} & 0 \\
               0 & 0 & -d\delta_{2}
             \end{array}
          \right] \\
   &&  \hspace{3.5cm}-
          \left [
             \begin{array}{ccc}
               \lambda_{1} & 0 & 0 \\
               0 & -\delta_{1} & 0 \\
               0 & 0 & -\delta_{2}
             \end{array}
          \right]
          \left [
             \begin{array}{ccc}
               0 & -h'_{2}dh_{1} & -h'_{3}dh_{1} \\
               h'_{2}dh_{1} & 0 & -h'_{3}dh_{2} \\
               h'_{3}dh_{1} & h'_{3}dh_{2} & 0
             \end{array}
          \right] \\
   &=& \left [
             \begin{array}{ccc}
               0 & \delta_{1} h'_{2}dh_{1} & \delta_{2} h'_{3}dh_{1} \\
               \lambda_{1} h'_{2}dh_{1} & 0 & \delta_{2} h'_{3}dh_{2} \\
               \lambda_{1} h'_{3}dh_{1} & -\delta_{1} h'_{3}dh_{2} & 0
             \end{array}
          \right]+
          \left [
             \begin{array}{ccc}
               d\lambda_{1} & 0 & 0 \\
               0 & -d\delta_{1} & 0 \\
               0 & 0 & -d\delta_{2}
             \end{array}
          \right] \\
   && \hspace{4cm}-
         \left [
             \begin{array}{ccc}
               0 & -\lambda_{1} h'_{2}dh_{1} & -\lambda_{1} h'_{3}dh_{1} \\
               -\delta_{1} h'_{2}dh_{1} & 0 & \delta_{1}h'_{3}dh_{2} \\
               -\delta_{2} h'_{3}dh_{1} & -\delta_{2} h'_{3}dh_{2} & 0
             \end{array}
          \right]\\
    &=& \left [
             \begin{array}{ccc}
               d\lambda_{1} & (\lambda_{1}+ \delta_{1})h'_{2}dh_{1} & (\lambda_{1}+\delta_{2}) h'_{3}dh_{1} \\
               (\lambda_{1}+ \delta_{1}) h'_{2}dh_{1} & -d\delta_{1} & (\delta_{2}-\delta_{1})h'_{3}dh_{2} \\
               (\lambda_{1}+ \delta_{2}) h'_{3}dh_{1} & (-\delta_{1}+\delta_{2}) h'_{3}dh_{2} & -d\delta_{2}
             \end{array}
          \right].
\end{eqnarray*}
We know that $(\mathbf{H}'d\mathbf{AH}) = (d\mathbf{A})$, then a column by column
computation of the exterior product of the subdiagonal elements of
$\mathbf{H}'d\mathbf{HD} + d\mathbf{D} - \mathbf{DH}'d\mathbf{H}$ gives, ignoring the
sign,
$$
  (d\mathbf{A}) = (\lambda_{1}+ \delta_{1})(\lambda_{1}+ \delta_{2})(-\delta_{1}+\delta_{2})\left (\bigwedge_{i =1}^{3}\bigwedge_{j =i +
   1}^{3} h'_{j}dh_{i} \right )\wedge d\lambda_{1}\wedge -d \delta_{1}\wedge
   -d\delta_{2}.
$$
Recall that, if for example, the first element in each column of $\mathbf{H}$ is
nonnegative, so, the transformation $\mathbf{A} = \mathbf{HDH}'$ is $1-1$. Then the
corresponding jacobian must be divided by $2^{m}$, see \citet[pp. 104-105]{m:82}.
Thus we have
$$
  (d\mathbf{A}) = 2^{-3}(\lambda_{1}+ \delta_{1})(\lambda_{1}+ \delta_{2})(\delta_{1}-\delta_{2})
         (\mathbf{H}'d\mathbf{H})\wedge (d\mathbf{D}),
$$
where $(\mathbf{H}'d\mathbf{H})$ is the Haar measure on $\mathcal{O}(m)$ and
$$
  (\mathbf{H}'d\mathbf{H}) = \bigwedge_{i < j}^{m} h'_{j}dh_{i}, \quad (d\mathbf{D}) =
  d\lambda_{1}\wedge d \delta_{1}\wedge
   d\delta_{2},
$$
this is, $(d\mathbf{D})$ is a exterior product of all differentials $d\lambda_{i}$
and $d\delta_{j}$ ignoring the sign.

Analogously, if $\mathbf{A} \in \mathcal{S}_{3}^{\pm}(2,1)$,
$$
  (d\mathbf{A}) = 2^{-3}(\lambda_{1}- \lambda_{2})(\lambda_{1}+ \delta_{1})(\lambda_{2}+\delta_{1})
         (\mathbf{H}'d\mathbf{H})\wedge (d\mathbf{D}).
$$
Similarly, let $\mathbf{A} \in \mathcal{S}_{4}^{\pm}(2,2)$, then
$$
  (d\mathbf{A}) = 2^{-4}(\lambda_{1}- \lambda_{2})(\delta_{1}- \delta_{2})(\lambda_{1}+
        \delta_{1})(\lambda_{1}+\delta_{2})(\lambda_{2}+ \delta_{1})(\lambda_{2}+
        \delta_{2})(\mathbf{H}'d\mathbf{H})\wedge (d\mathbf{D}).
$$
By mathematical induction we have
\begin{thm}
Let $\mathbf{A} \in \mathcal{S}_{m}^{\pm}(m_{1}, m_{2})$ such that $\mathbf{A} =
\mathbf{HDH}'$, where $\mathbf{H} \in \mathcal{O}(m)$, $\mathbf{D}$ is a diagonal
matrix with $\lambda_{1}>\cdots> \lambda_{m_{1}}> 0$ and $0 > -\delta_{1}>\cdots>
-\delta_{m_{2}}$, $m_{1}+m_{2} = m$. Then
$$
  (d\mathbf{A}) = 2^{-m} \prod_{i < j}^{m_{1}}(\lambda_{i}- \lambda_{j})\prod_{i < 1}^{m_{1}}(\delta_{i}- \delta_{j})
    \prod_{i,j}^{m_{1}, m_{2}}(\lambda_{i}+ \delta_{j})(\mathbf{H}'d\mathbf{H})\wedge (d\mathbf{D}).
$$
where
$$
   \prod_{i,j}^{m_{1}, m_{2}}(\lambda_{i}+ \delta_{j}) =  \prod_{i=1}^{m_{1}}
    \prod_{j=1}^{m_{2}}(\lambda_{i}+ \delta_{j}), \quad (\mathbf{H}'d\mathbf{H})
    = \bigwedge_{i < j}^{m} h'_{j}dh_{i}, \quad (d\mathbf{D})= \bigwedge_{i=1}^{m_{1}} d\lambda_{i}
    \bigwedge_{j=1}^{m_{2}} d\delta_{j}.
$$
\end{thm}
A similar procedure for $\mathbf{A} \in \mathcal{S}_{m}^{\pm}(q,q_{1}, q_{2})$ gives:
\begin{thm}
Let $\mathbf{A} \in \mathcal{S}_{m}^{\pm}(q,q_{1}, q_{2})$ such that $\mathbf{A} =
\mathbf{H}_{1}\mathbf{D}\mathbf{H}'_{1}$, where $\mathbf{H}_{1} \in
\mathcal{V}_{q,m}$, $\mathbf{D}$ is a diagonal matrix with $\lambda_{1}>\cdots>
\lambda_{q_{1}}> 0$ and $0 > -\delta_{1}>\cdots> -\delta_{q_{2}}$, $q_{1}+q_{2} =q$.
Then
$$
  (d\mathbf{A}) = 2^{-q}\prod_{i = 1}^{q_{1}}\lambda_{i}^{m-q}  \prod_{j = 1}^{q_{2}}\delta_{j}^{m-q}
    \prod_{i < j}^{q_{1}}(\lambda_{i}- \lambda_{j})\prod_{i < 1}^{q_{1}}(\delta_{i}- \delta_{j})
    \prod_{i,j}^{q_{1}, q_{2}}(\lambda_{i}+ \delta_{j})(\mathbf{H}'d\mathbf{H})\wedge (d\mathbf{D}).
$$
where
$$
   \prod_{i,j}^{q_{1}, q_{2}}(\lambda_{i}+ \delta_{j}) =  \prod_{i=1}^{q_{1}}
    \prod_{j=1}^{q_{2}}(\lambda_{i}+ \delta_{j}), \quad (\mathbf{H}'_{1}d\mathbf{H}_{1}) = \bigwedge_{i=1}^{m}
    \bigwedge_{j = i+1}^{q}h'_{j}dh_{i}, \quad (d\mathbf{D})= \bigwedge_{i=1}^{q_{1}} d\lambda_{i}
    \bigwedge_{j=1}^{q_{2}} d\delta_{j}.
$$
\end{thm}

\section{Jacobians of symmetric matrices with multiplicity in its eigenvalues} \label{multi}

As a motivation of this section, consider a general random matrix $\mathbf{A} \in
\mathfrak{R}^{m \times m}$, explicitly
$$
  \mathbf{A} = \left [
         \begin{array}{ccc}
           a_{11} & \cdots & a_{1m} \\
           \vdots & \ddots & \vdots \\
           a_{m1} & \cdots & a_{mm}
         \end{array}
      \right].
$$
Any density function of this matrix can be expressed as
$$
  dF_{\mathbf{A}}(\mathbf{A}) = f_{\mathbf{A}}(\mathbf{A})(d\mathbf{A}),
$$
where $(d\mathbf{A})$ denotes the measure of Lebesgue in $\mathfrak{R}^{m^{2}}$,
which can be written by using the exterior product, as
$$
  (d\mathbf{A}) = \bigwedge_{i=1}^{m} \bigwedge_{j=1}^{m} da_{ij},
$$
see \citet{m:82}.

However, if $\mathbf{A} \in {\mathcal S}_{m}$ and it is non singular,
  then the measure of Lebesgue defined in ${\mathcal S}_{m}$ is given by
\begin{equation}\label{sim}
    (d\mathbf{A}) = \bigwedge_{i\leq j}^{m} da_{ij}.
\end{equation}
\begin{rem}
Note that the above product is the measure of Hausdorff  on $\mathfrak{R}^{m^{2}}$
defined on the homogeneous space of positive definite symmetric matrices, see
\citet{b:86}
\end{rem}

In general, we can consider any factorization of the Lebesgue measure $(d\mathbf{A})$
on ${\mathcal S}_{m}$ as an alternative definition of $(d\mathbf{A})$ with respect to
the corresponding coordinate system. For example, if we consider the spectral
decomposition (SD), $\mathbf{A} = \mathbf{HDH}'$, where $\mathbf{H} \in
\mathcal{O}(m)$, and $\mathbf{D}$ is a diagonal matrix with $D_{1}>\cdots> D_{m}> 0$
or we consider the Cholesky decomposition $\mathbf{A} = \mathbf{T}'\mathbf{T}$, where
$\mathbf{T}$ is upper-triangular with positive diagonal elements, then we have
respectively
$$
  (d\mathbf{A}) = \left\{
           \begin{array}{ll}
             2^{-m}\Prod{i<j}{m}(D_{i}-D_{j})(\mathbf{H}'d\mathbf{H})\wedge (d\mathbf{D}),
                & \hbox{Spectral decomposition;} \\
             2^{m}\Prod{i=1}{m} t_{ii}^{m+1-i}(d\mathbf{T}), & \hbox{Cholesky decomposition,}
           \end{array}
         \right .
$$
see \citet{dgggj:04a}, where
$$
  (\mathbf{H}'d\mathbf{H}) = \bigwedge_{i < j}^{m} h'_{j}dh_{i}, \quad (d\mathbf{D})= \bigwedge_{i
  =1}^{m} dD_{i} \mbox{\ and \ } (d\mathbf{T}) = \bigwedge_{i \leq j}^{m} dt_{ij}.
$$

In some occasions is difficult to establish an explicit form of the Lebesgue o
Hausdorff measures in the original coordinate system. In particular if $\mathbf{A}
\in \mathcal{S}_{m}^{+}(q)$ some unsuccessful efforts have been trailed, see
\citet{s:03} and \citet{dg:07}. A definition of such measure in terms of the SD is
given by \citet{u:94}:
\begin{equation}\label{SD}
    (d\mathbf{A}) = 2^{-q} \prod_{i=1}^{q} D_{i}^{m-q} \prod_{i<j}^{q} (D_{i}-D_{j})
        (\mathbf{H}_{1}d\mathbf{H}_{1})\wedge (d\mathbf{D}),
\end{equation}
where $\mathbf{H}_{1} \in \mathcal{V}_{q,m}$, $\mathbf{D}$ is a diagonal matrix with
$D_{1}>\cdots> D_{q}> 0$ and
$$
  (\mathbf{H}'_{1}d\mathbf{H}_{1}) = \bigwedge_{i =1}^{m}\bigwedge_{j =i+1}^{q} h'_{j}dh_{i},
    \quad (d\mathbf{D})= \bigwedge_{i=1}^{q} dD_{i} ;
$$
for alternative expressions of $(d\mathbf{A})$ in terms of other factorizations see
\citet{dgggj:04a} and \citet{dgggd:04b}.

Now suppose that one (or more) eigenvalue(s) of $\mathbf{A} \in \mathcal{S}_{m}$ has
(have) multiplicity. Then consider $\mathbf{A} = \mathbf{HDH}'$, where $\mathbf{H}
\in \mathcal{O}(m)$, $\mathbf{D}$ is a diagonal matrix with $D_{1}\geq \cdots \geq
D_{m}> 0$. Moreover, let $D_{k_{1}}, \dots D_{k_{l}}$ be the $l$ distinct eigenvalues
of $A$, i.e. $D_{k_{1}}> \cdots > D_{k_{l}}
> 0$, where $m_{j}$ denotes the repetitions of the eigenvalue $D_{k_{j}}$, $j
=1, 2, \dots, l$, and of course $m_{1} + \cdots m_{l} = m$; finally denote the
corresponding set of matrices by $\mathbf{A} \in \mathcal{S}_{m}^{l}$. It is clear
that $\mathbf{A}$ exists in the homogeneous subspace of the symmetric matrices with
dimension $m(m+1)/2$ of rank $m$; more accurately, when there exist multiplicity in
the eigenvalues, $\mathbf{A}$ exists in the manifold of dimension $ml-l(l-1)/2$, even
exactly for computations we say that $\mathbf{A} \in \mathcal{S}_{m}^{+}(l)$. For
proving it, consider the matrix $\mathbf{B} \in \mathcal{S}_{2}$, such that
$\mathbf{S} = \mathbf{HDH}'$, here $\mathbf{H} \in \mathcal{O}(2)$, and $\mathbf{D}$
is a diagonal matrix with $D_{1}\geq D_{2} > 0$ where $D_{1}= D_{2} = \kappa$, then
the measure
$$
  (d\mathbf{B}) = 2^{-2}\Prod{i<j}{2}(D_{i}-D_{j})(\mathbf{H}'d\mathbf{H})\wedge d\mathbf{D} =
2^{-2}(\kappa-\kappa)(\mathbf{H}'d\mathbf{H})\wedge d\mathbf{D} = 0.
$$
Also note that, in fact the measure $(d\mathbf{D}) = dD_{1}\wedge dD_{2} =
d\kappa\wedge d\kappa = 0$.  This is analogous to the following situation, to propose
for a curve in the space $(\mathfrak{R}^{3})$ the measure of Lebesgue defined by
$dx_{1}\wedge dx_{2}$.

Now, when we consider the factorization of the measure of Lebesgue in terms of the
spectral decomposition, we do not have $2(2+1)/2 = 3$ but only $2(1)-1(1+1)/2 + 1 =
2$ mathematical independent elements in $\mathbf{B}$, because in $\mathbf{D}$,
$D_{1}= D_{2} = \kappa$ and then there is only one mathematical independent element.

Also, observe that the space of positive definite $m \times m$ matrices is a subset
of Euclidian space of symmetric $m \times m$ matrices of dimension $m(m+1)/2$, and in
fact it forms an open cone described by the following system of inequalities, see
\citet[p. 61 and p. 77 Problem 2.6]{m:82}:
\begin{equation}\label{cono}
    \mathbf{A} > 0 \Leftrightarrow a_{11}>0, \det \left[ \begin{array}{cc}
                                                a_{11} & a_{12} \\
                                                a_{21} & a_{22}
                                              \end{array}
    \right] > 0, \cdots, \det (\mathbf{A}) > 0.
\end{equation}
In particular, let $m = 2$, after factorizing the measure of Lebesgue in
$\mathcal{S}_{m}$ by the spectral decomposition, the inequalities (\ref{cono}) are as
follows
\begin{equation}\label{cono2}
    \mathbf{A} > 0 \Leftrightarrow D_{1} >0, D_{2} >0, D_{1}D_{2}>0.
\end{equation}
But if  $D_{1}= D_{2} = \kappa$, (\ref{cono2}) it reduces to
\begin{equation}\label{cono3}
    \mathbf{A} > 0 \Leftrightarrow \kappa >0, \kappa^{2}>0.
\end{equation}
Which defines a curve (a parabola) in the space, over the line $D_{1}= D_{2} (=
\kappa)$ in the subspace of points $(D_{1},D_{2})$.

A similar situation appear in the following cases: i) When we consider multiplicity
of the singular values in the SVD; such set of matrices will be denoted by
$\mathbf{X} \in {\mathcal L}_{m,N}^{l}(q)$, $q \geq l$; or by $\mathbf{X} \in
{\mathcal L}_{m,N}^{+}(q,l)$  $q \geq l$;  ii) If we consider multiplicity in the
eigenvalues; the corresponding set of matrices will be denoted by $\mathbf{A} \in
\mathcal{S}_{m}^{+}(q,l)$, $q \geq l$; iii) And if $A$ is nonpositive definite.

Thus, unfortunately, we must qualify as incorrect the asseverations of \citet{z:07}
about the validity of his Lemmas 2, 3 and  consequences, under multiplicity
assumptions of singular values and eigenvalues.

As a summary we have the next results, which collect the main conclusions of
 Section \ref{indef} and the present section, and follow a similar proof of Theorem 1 in
\citet{dgg:06}:
\begin{thm}
Consider $\mathbf{Y} \in {\mathcal L}_{m,N}(q)$ and  $\mathbf{Y} = \mathbf{X}^{+}$, then%
$$
  (d\mathbf{Y}) = \prod_{i = 1}^{k}\sigma_{i}^{-2(N+m-k)}(d\mathbf{X})
$$
where $\mathbf{X} = \mathbf{H}_{1}\mathbf{D}_{\sigma} \mathbf{P}'_{1}$ is the
nonsingular part of SVD of $\mathbf{X} $, with $\mathbf{H}_{1} \in
\mathcal{V}_{k,N}$, $\mathbf{P}_{1} \in \mathcal{V}_{k,m}$, $\mathbf{D}_{\sigma} =
\diag(\sigma_{1}, \dots, \sigma_{k})$, $\sigma_{1}> \cdots > \sigma_{k}
> 0$, the measure $(d\mathbf{X})$ is
$$
  (d\mathbf{X}) = 2^{-k} \prod_{i = 1}^{k}\sigma_{i}^{(N+m-2k)}\prod_{i <j}^{k}(\sigma_{i}^{2}- \sigma_{j}^{2})
        (\mathbf{H}'_{1}d\mathbf{H}_{1})\wedge(\mathbf{P}'_{1}d\mathbf{P}_{1})\wedge (d\mathbf{D}_{\sigma}),
$$
and
$$
  k =
  \left\{%
     \begin{array}{ll}
        q, & \hbox{$\mathbf{X} \in {\mathcal L}_{m,N}^{+}(q)$;} \\
        l, & \hbox{$\mathbf{X} \in {\mathcal L}_{m,N}^{+}(q,l)$.} \\
    \end{array}%
  \right.
$$
\end{thm}

Similarly, for symmetric matrices we have,

\begin{thm}
Let $\mathbf{V} \in \mathfrak{R}^{m \times m}$ be a symmetric matrix and let
$\mathbf{W} = \mathbf{V}^{+}$, then
\begin{enumerate}
\item
$$
  (d\mathbf{W}) = \prod_{i = 1}^{\beta}|\lambda_{i}|^{-2m + \beta - 1}(d\mathbf{V}),
$$
where $\mathbf{V} = \mathbf{H}_{1}\mathbf{D}_{\lambda} \mathbf{H}'_{1}$ is the
nonsingular part of SD of $\mathbf{V}$, with $\mathbf{H}_{1} \in
\mathcal{V}_{\beta,N}$, $\mathbf{D}_{\lambda} = \diag(\lambda_{1}, \dots,
\lambda_{\beta})$, $|\lambda_{1}|> \cdots > |\lambda_{\beta}| > 0$, the measure
$(d\mathbf{V})$ is
$$
  (d\mathbf{V}) = 2^{-\beta} \prod_{i=1}^{\beta} |\lambda_{i}|^{m-\beta} \prod_{i<j}^{\beta}
    (|\lambda_{i}|-|\lambda_{j}|)(\mathbf{H}_{1}d\mathbf{H}_{1})\wedge
    (d\mathbf{D}_{\lambda}),
$$
and
$$
  \beta =
  \left\{%
     \begin{array}{ll}
        m, & \hbox{$\mathbf{V}$ or $-\mathbf{V} \in \mathcal{S}_{m}$;} \\
        l, & \hbox{$\mathbf{V}$ or $-\mathbf{V} \in \mathcal{S}_{m}^{l}$;} \\
        q, & \hbox{$\mathbf{V}$ or $-\mathbf{V} \in \mathcal{S}_{m}^{+}(q)$;} \\
        k, & \hbox{$\mathbf{V}$ or $-\mathbf{V} \in \mathcal{S}_{m}^{+}(q,k)$.} \\
    \end{array}%
  \right.
$$
\item
$$
  (d\mathbf{W}) = \prod_{i = 1}^{\alpha_{1}}\lambda_{i}^{-2(m - \alpha_{1}/2 - \alpha_{2} + 1)}
             \prod_{j = 1}^{\alpha_{2}}\delta_{j}^{-2(m - (\alpha - 1)/2)}(d\mathbf{V}),
$$
where $\alpha = \alpha_{1} + \alpha_{2}$, $\mathbf{V} = \mathbf{H}_{1}\mathbf{D}
\mathbf{H}'_{1}$ is the nonsingular part of SD of $\mathbf{V}$, with $\mathbf{H}_{1}
\in \mathcal{V}_{\alpha,N}$, $\mathbf{D} = \diag(\lambda_{1}, \dots,
\lambda_{\alpha_{1}}, -\delta_{1}, \dots, -\delta_{\alpha_{1}})$, $\lambda_{1}>
\cdots > \lambda_{\alpha_{1}}
> 0$; $|\delta_{1}|> \cdots > |\delta_{\alpha_{2}}| > 0$, the measure $(d\mathbf{V})$ is
$$
  (d\mathbf{V}) = 2^{-\alpha}\prod_{i = 1}^{\alpha_{1}}\lambda_{i}^{m-\alpha}
    \prod_{j = 1}^{\alpha_{2}}\delta_{j}^{m-\alpha}\prod_{i < j}^{\alpha_{1}}
    (\lambda_{i}- \lambda_{j})\prod_{i < 1}^{\alpha_{1}}(\delta_{i}- \delta_{j})
    \prod_{i,j}^{\alpha_{1}, \alpha_{2}}(\lambda_{i}+ \delta_{j})(\mathbf{H}'
    d\mathbf{H})\wedge (d\mathbf{D}).
$$
and
$$
   \alpha =
  \left\{%
     \begin{array}{ll}
        m, & \hbox{$\mathbf{V} \in \mathcal{S}^{\pm}_{m}(m_{1},m_{2})$;} \\
        l, & \hbox{$\mathbf{V} \in \mathcal{S}^{\pm}_{m}(l_{1},l_{2})$;} \\
        q, & \hbox{$\mathbf{V} \in \mathcal{S}^{\pm}_{m}(q, q_{1},q_{2})$;} \\
        k, & \hbox{$\mathbf{V} \in \mathcal{S}^{\pm}_{m}(q, k_{1},k_{2})$,} \\
    \end{array}%
  \right.
$$
and $\mathbf{V} \in \mathcal{S}^{\pm}_{m}(l_{1},l_{2})$ denotes a nonsingular
indefinite matrix with multiplicity in its eigenvalues and $\mathbf{V} \in
\mathcal{S}^{\pm}_{m}(q, k_{1},k_{2})$ denotes a singular indefinite matrix with
multiplicity in its eigenvalues.
\end{enumerate}
\end{thm}

\section{Conclusions}

This work determines the jacobians of the SVD and the SD under multiplicity of the
singular values and eigenvalues, respectively. For the SD case, we compute the
jacobian for nonsingular and singular indefinite matrices with and without
multiplicity in their eigenvalues. Also, we calculate the jacobians for a general
matrix and its Moore-Penrose inverse, and for a symmetric matrix with all its
variants (nonpositive, nonnegative and indefinite). In every case we specify the
measures of Hausdorff which support the jacobian computations. These results detecte
and correct some inconsistences in the validity of Lemmas 2 and 3 by \citet{z:07},
this is, the proofs in \citet{z:07} for Lemmas 2 and 3 are valid only when all the
nonzero singular values or eigenvalues are distinct and for positive or negative
semidefinite matrices, as it is assumed in \citet{dgt:05} and \citet{dgg:06}. We
highlight that the results of this paper will be the foundations of an explored
problem in literature: the test criteria in MANOVA when there exist multiplicities
in: the matrix of sum of squares and sum of products, due to the hypothesis
$\mathbf{S}_{H}$; the matrix of sum of square and sum of products, due to the error
$\mathbf{S}_{E}$; the matrices $\mathbf{S}_{H}\mathbf{S}_{E}^{-1}$; $(\mathbf{S}_{H}
+ \mathbf{S}_{E})^{-1}\mathbf{S}_{H}$; see \citet{dgc:08}.

\section*{Acknowledgment}

This research work was partially supported  by IDI-Spain, grants FQM2006-2271 and
MTM2008-05785, and CONACYT-M\'exico, research grant no. \ 81512. This paper was
written during J. A. D\'{\i}az- Garc\'{\i}a's stay as a visiting professor at the
Department of Statistics and O. R. of the University of Granada, Spain.


\end{document}